\documentclass{amsart}

\usepackage{amssymb}
\newtheorem{theorem}{Theorem}[section]
\newtheorem{lemma}[theorem]{Lemma}
\newtheorem*{conjecture}{Conjecture}

\theoremstyle{definition}

\theoremstyle{remark}

\numberwithin{equation}{section}

\begin{document}

\title[RANDOM HARMONIC POLYNOMIALS]{A note on the expectation of zeros of random harmonic polynomials: The Kac model}
\thanks{This is the accepted author manuscript. First published in
\emph{Proceedings of the American Mathematical Society} in volume 153,
number 2 (2025), published by the American Mathematical Society.
Pages 611--623. DOI: \texttt{10.1090/proc/17002}.
\textcopyright\ 2024 American Mathematical Society.}


\author{Dawei Lu}
\address{School of Mathematical Sciences, Dalian University of Technology, Dalian 116023, China}
\email{ludawei\_dlut@163.com}
\thanks{This first author was supported by National Key R\&D Program of China (No.2023YFA1009200), National Natural Science Foundation of China (Grant No.12371154), the Fundamental Research Funds for the Central Universities (DUT23LAB307), Dalian High-level Talent Innovation Program (Grant No.2020RD09).}

\author{Yuchen Wang}
\address{School of Mathematical Sciences, Dalian University of Technology, Dalian 116023, China}
\email{wycwyc@mail.dlut.edu.cn}

\subjclass[2020]{Primary 30C15, 34F05, 60G15; Secondary 26C10}

\date{}

\begin{abstract}
Motivated by  the questions posed by W. V. Li and A. Wei and the conjecture of E. Lundberg and A. Thomack, we study the expected number of zeros of random harmonic polynomials $H_{n,m}(z)= p_{n}(z)+\overline{q_{m}(z)}$ with independently and identically distributed Gaussian coefficients. In this paper we verify the conjecture of E. Lundberg and A. Thomack that the expectation is $O(n)$ when $\deg p = \alpha \deg q$, where $0\leq\alpha<1$. This result extends the previous estimates when $m$ is a fixed constant or $m=n$ to more general case.
\end{abstract}

\maketitle


\bibliographystyle{amsplain}
\section{Introduction}
A random harmonic polynomial is a complex-valued function in $\mathbb{C}$ (the complex plane) of the form $H_{n,m}(z)= p_{n}(z)+\overline{q_{m}(z)}$, where $p_{n}(z)$ and $q_{m}(z)$ are both analytic polynomials, with $n=\deg p\geq \deg q =m$. Since $H_{n,m}(z)$ is not an analytic function, the Fundamental Theorem of Algebra is not applicable. Consequently, it becomes pertinent to inquire about the number of zeros that $H_{n,m}(z)$ can possess. In 1992, T. Sheil-Small initiated the investigation into the maximum number of zeros a polynomial can have and hypothesized that the maximum number of zeros for $H_{n,m}(z)$ is no larger than $n^{2}$. A. S. Wilmshurst proved this in \cite{Wilmshurst1994}, \cite{wilmshurst1998valence} using Bézout's Theorem that the maximum is at most $n^{2}$. He also illustrated with examples that the upper bound is sharp for cases where $m=n$ or $m=n-1$. In \cite{wilmshurst1998valence}, A. S. Wilmshurst conjectured that the maximum number of zeros is $m(m-1)+3n-2$ for $m\leq n-1$. D. Khavinson and G. \' Swiatek \cite{Khavinson2002} successfully proved that the maximum is $3n-2$ when $m=n-1$, utilizing an indirect approach grounded in Fatou's Theorem from the field of holomorphic dynamics. And L. Geyer \cite{geyer2008sharp} proved the sharpness of the upper bound. Other related works can be found in \cite{bshouty1995exact}, \cite{bshouty2004crofoot}, \cite{khavinson2024valence}, \cite{khavinson2006number}, \cite{lee2017new}, and \cite{lundberg2023valence}. Considering the fluctuating count of zeros for $H_{n,m}(z)$, which is different from polynomials of fixed coefficients, it is natural to understand average case behaviour on $\mathbb{C}$ with random coefficients in addition to extremal behaviour.

There is a long history of studying zeros of a random polynomial whose coefficients are independent. These early works include \cite{Bloch1932} and \cite{kac1943average}. More details are provided in \cite{edelman1995many}. Exact formula for the expected number of real zeros under independent identically distributed Gaussian coefficients are found for a random polynomial by Kac \cite{kac1943average}, and for a random trigonometric polynomial by Dunnage \cite{dunnage1966number}. In \cite{shepp1995complex}, \cite{ibragimov1997roots} and \cite{Peres2005},  there is also a significant amount of research presented on the complex roots over a fixed domain. Especially pertinent to this paper is our observation that the study on the probability of harmonic polynomials, started by W. V. Li and A. Wei \cite{Li2009}, which is similar to studies on stochastic gravitational lensing \cite{bleher2014counting}, \cite{khavinson2010transcendental}, \cite{petters2009mathematical}, \cite{petters2009mathematical2} and \cite{wei2017numbers}. W. V. Li and A. Wei \cite{Li2009} utilized similar methodologies, including the Kac-Rice formula applied to vector fields. They considered the harmonic Kostlan ensemble where the harmonic polynomial $H(z)=H_{n,m}(z)= p_{n}(z)+\overline{q_{m}(z)}$ is randomized by sampling $p_{n},q_{m}$ independently from the analytic Kostlan ensemble. In this situation, they successfully reformulated the Kac-Rice formula, resulting in a more convenient computation process. Moreover, they also provided the formula for calculation in the case of independent and identically distributed Gaussian coefficients.

\begin{theorem}[Theorem 4.1 in W. V. Li and A. Wei \cite{Li2009}]
The expectation \( \mathbb{E}[N_{H}(T)] \) of the numbers of zeros in an open set \( T \subseteq \mathbb{C} \) of a harmonic polynomial \( H(z) = p_{n}(z) + \overline{q_{m}(z)} \)  with i.i.d. Gaussian coefficients satisfies
\[ \mathbb{E}[N_{H}(T)] = \frac{1}{\pi} \int_{T} \frac{1}{\vert z\vert^{2}} \frac{r_{1}^{2}+r_{2}^{2}-2r_{12}^{2}}{r_{3}^{2}\sqrt{(r_{1}+r_{2})^{2}-4r_{12}^{2}}} \, dA(z), \]
where \( dA(z) \) denotes the Lebesgue measure on the plane, and
\begin{flalign}
    &r_{3}=\sum_{j=0}^{n}\vert z \vert^{2j}+\sum_{j=0}^{m}\vert z \vert^{2j},  
    &r_{12}=\left(\sum_{j=0}^{n}j\vert z \vert^{2j}\right)\left(\sum_{j=0}^{m}j\vert z \vert^{2j}\right),  \notag
    && \\
    &r_{1}=r_{3}\sum_{j=0}^{n}j^{2}\vert z \vert^{2j}-\left(\sum_{j=0}^{n}j\vert z \vert^{2j}\right)^{2},  
    &r_{2}=r_{3}\sum_{j=0}^{m}j^{2}\vert z \vert^{2j}-\left(\sum_{j=0}^{m}j\vert z \vert^{2j}\right)^{2}.  \notag
    &&
\end{flalign}
\end{theorem}

This result is proved using a version of the Kac-Rice formula for vector fields which gives
\[ \mathbb{E}[N_{H}(\mathbb{T})] = \int_{T} \mathbb{E}(\lvert \det J_{H}(z) \rvert \big\vert H(z) = 0) \rho(0; z) \, dA(z), \]
where, for each \(z\), \(\rho(u; z)\) is the probability density function of \(u = H(z)\).\\

Additionally, their result is further proved to be valid for all random harmonic polynomials with Gaussian coefficients. For more details, please refer to \cite{Lerario2016}, \cite{Thomack2018} and \cite{Lundberg2023}. Especially, E. Lundberg and A. Thomack obtained the following theorem of a general case in \cite{Lundberg2023}.
\begin{lemma}[Theorem 2.1 in E. Lundberg and A. Thomack \cite{Lundberg2023}]
\label{lemma1}
The expectation $\mathbb{E}[N_{H}(T)]$ of the number of zeros of
\[H_{n,m}(z)=\sum_{j=0}^{n}A_{j}z^{j}+\sum_{j=0}^{m}B_{j}\overline{z}^{j}\]
where $A_{0},\dots,A_{n}$ and $B_{0},\dots,B_{m}$ are mutually independent complex Gaussian random variables with $\mathbb{E}A_{j}=\mathbb{E}B_{j}=\mathbb{E}[A_{j}^{2}]=\mathbb{E}[B_{j}^{2}]=0$ and $\mathbb{E}A_{j}\overline{A_{j}}=\alpha_{j}$ and $\mathbb{E}B_{j}\overline{B_{j}}=\beta_{j}$ on a domain $T\subset \mathbb{C}$ is given by:
\[ \mathbb{E}[N_{H}(\mathbb{T})] = \frac{1}{\pi} \int_{T} \frac{1}{\vert z\vert^{2}} \frac{r_{1}^{2}+r_{2}^{2}-2r_{12}^{2}}{r_{3}^{2}\sqrt{(r_{1}+r_{2})^{2}-4r_{12}^{2}}} \, dA(z), \]
where \( dA(z) \) denotes the Lebesgue measure on the plane, and
\begin{flalign}
    &r_{3}=\sum_{j=0}^{n}\alpha_{j}\vert z \vert^{2j}+\sum_{j=0}^{m}\beta_{j}\vert z \vert^{2j},  
    &r_{12}=\left(\sum_{j=0}^{n}j\alpha_{j}\vert z \vert^{2j}\right)\left(\sum_{j=0}^{m}j\beta_{j}\vert z \vert^{2j}\right),  \notag
    && \\
    &r_{1}=r_{3}\sum_{j=0}^{n}j^{2}\alpha_{j}\vert z \vert^{2j}-\left(\sum_{j=0}^{n}j\alpha_{j}\vert z \vert^{2j}\right)^{2},  
    &r_{2}=r_{3}\sum_{j=0}^{m}j^{2}\beta_{j}\vert z \vert^{2j}-\left(\sum_{j=0}^{m}j\beta_{j}\vert z \vert^{2j}\right)^{2}.  \notag
    &&
\end{flalign}
\end{lemma}

For some early results in \cite{Lerario2016}-\cite{Li2009} and \cite{Thomack2016}-\cite{Thomack2018}, by applying the aforementioned lemma, the expectation of zeros of random harmonic polynomials in numerous distinct scenarios can be re-derived. In this paper, we sample $p_{n}$, $q_{m}$ independently from the analytic Kac ensemble, i.e. $A_{j},B_{j}\sim N_{\mathbb{C}}(0,1)$ are independent identically distributed standard complex Gaussian variables. Regarding this case, W. V. Li and A. Wei \cite{Li2009} conjectured that the average number of zeros is asymptotically $n$ when $m$ is fixed as $n\rightarrow \infty$. This was proved by A. Thomack in \cite{Thomack2016} where it was further conjectured that when $m\sim n$ as $n\rightarrow \infty$ the expectation of zeros satisfies $\mathbb{E}N_{H}(\mathbb{C})\sim Cn\log n$ for some constant $C$. In \cite{Lundberg2023}, E. Lundberg and A. Thomack proved that $\mathbb{E}N_{H}(\mathbb{C})\sim \frac{1}{2}n\log n$ when $m= n$ as $n\rightarrow \infty$ and provided the following conjecture.

\begin{conjecture}
In the case of independent and identically distributed Gaussian coefficients and $m=\alpha n$, where $0<\alpha<1$, the expected number of zeros of \( H(z) = p_n(z) + \overline{q_m(z)} \) is asymptotically proportional to \( n \), with the constant of proportionality perhaps depending on \( m \).
\end{conjecture}

In the process of proving this conjecture, we encountered the following difficulties. Firstly, when the coefficients are independent and identically distributed, compared to the case of $m $ is a fixed constant or $m=n$, we need to estimate more terms in the case of $m=\alpha n$, where $0\leq\alpha<1$. Thereby the calculations become more complex. Secondly, the series $\sum_{j=0}^{n}\vert z \vert^{2j}$ and $\sum_{j=0}^{m}\vert z \vert^{2j}$ force us to be more cautious when dealing with estimates near $\vert z\vert=1$ to avoid the use of $\sum_{i=0}^{n}x^{i}=\frac{x^{n+1}-1}{x-1}$. Lastly, different from the Kostlan or Weyl ensembles, for the Kac ensemble, when $\vert z \vert$ is sufficiently large, the influence of $q_{m}$ becomes significantly smaller compared to $p_{n}$. Therefore, we need a redivision of the entire complex plane and achieve more precise estimates when $\vert z \vert$ approaches 1.

The remaining sections are organized as follows: The main results are presented in Section 2. In Section 3, we first conduct some preliminary calculations, and then provide the proofs of the main results.
\section{Main results}
Based on Lemma \ref{lemma1}, we estimate the expected number of zeros of $H(z) = p_n(z) + \overline{q_m(z)}$ , where \( p_n \) and \( q_m \) are independently sampled from the Kac ensemble.
\begin{theorem}
\label{theorem}
Let \( H(z) = p_{n}(z) + \overline{q_{m}(z)} \) be a random harmonic polynomial with independent and identically distributed Gaussian coefficients, where \( m = \alpha n ,  \alpha \in (0, 1) \). As \( n \rightarrow \infty \), the expectation of the number of zeros of \( H \) satisfies
\[ \mathbb{E}[N_{H}(\mathbb{C})] = Cn, \]
where C is a constant related to $\alpha$.
\end{theorem}
This theorem not only resolves the conjecture of E. Lundberg and A. Thomack, but also offers further insights into the expectation of zeros across different areas of the complex plane. Essentially, Theorem \ref{theorem} stems from upper bounds in independent regions and the lower bound in the entire complex plane as follow.
\begin{align}
\label{1}
\mathbb{E}[N_{H}(\lbrace \lvert z \rvert^{2} \leq 1-\frac{1}{n} \rbrace)] \,\lesssim\, n,
\end{align}
\begin{align}
\label{2}
\mathbb{E}[N_{H}(\lbrace 1-\frac{1}{n} \leq \lvert z \rvert^{2} \leq 1+\frac{1}{n} \rbrace)] \,\lesssim\, n , 
\end{align}
\begin{align}
\label{3}
\mathbb{E}[N_{H}(\lbrace \lvert z \rvert^{2} \geq 1+\frac{1}{n} \rbrace)] \,\lesssim\, n ,
\end{align}
\begin{align}
\label{4}
\mathbb{E}[N_{H}(\mathbb{C})] \,\geq\, n ,
\end{align}
where by $A\lesssim B$ we mean that there is a finite positive constant $K$ such that $A\leq KB$.

\section{Proofs}
Before giving the proofs of the main results, we first transform the Kac-Rice formula and make a simple scaling of the integrand. Applying Lemma 2.1 and integrating in polar coordinates $z=re^{i\theta}$ while making the change of variables $w=r^{2}$, $dw=2rdr$, we obtain the following identity for the expected number of zeros in an annulus $\Omega=\lbrace a < \lvert z \rvert^2 < b \rbrace$.
\begin{align*}
\mathbb{E}[N_{H}(\Omega)] &= \int_{a}^{b} \frac{1}{w} \frac{r_{1}^{2}+r_{2}^{2}-2r_{12}^{2}}{r_{3}^{2}\sqrt{(r_{1}+r_{2})^{2}-4r_{12}^{2}}} \, dw:= \int_{a}^{b} F(w) \, dw,
\end{align*}
where $r_{1}=(a_{m}+a_{n})c_{n}-b_{n}^{2}$, $r_{2}=(a_{m}+a_{n})c_{m}-b_{m}^{2}$, $r_{12}=b_{n}b_{m}$ , and $r_{3}= a_{n}+a_{m}$ with
\begin{align}
\label{abc}
a_{k}=\sum_{j=0}^{k}w^{j}, \quad b_{k}=\sum_{j=0}^{k}jw^{j}, \quad c_{k}=\sum_{j=0}^{k}j^{2}w^{j}.
\end{align}
Furthermore, we have
\begin{align}
\label{F}
F(w)= \frac{r_{1}^{2}+r_{2}^{2}-2r_{12}^{2}}{wr_{3}^{2}\sqrt{(r_{1}+r_{2})^{2}-4r_{12}^{2}}}\leq\frac{r_{1}-r_{2}}{wr_{3}^{2}}+\frac{\sqrt{r_{1}r_{2}-r_{12}^{2}}}{wr_{3}^{2}}.
\end{align} 
\subsection{Pre-calculation}
To facilitate subsequent proofs, we perform some preliminary calculations. First, by \eqref{abc}, we have
\begin{align}
\label{11}
\frac{b_{k}}{a_{k}}=\frac{\sum_{j=0}^{k}jw^{j}}{\sum_{j=0}^{k}w^{j}}\leq \frac{\sum_{j=0}^{k}kw^{j}}{\sum_{j=0}^{k}w^{j}}=k
\end{align}
and
\begin{align}
\label{12}
\frac{c_{k}}{a_{k}}=\frac{\sum_{j=0}^{k}j^{2}w^{j}}{\sum_{j=0}^{k}w^{j}}\leq \frac{\sum_{j=0}^{k}k^{2}w^{j}}{\sum_{j=0}^{k}w^{j}}=k^{2}.
\end{align}
Then, still utilizing \eqref{abc}, we have
\begin{align}
\label{13}
b_{k}=w\frac{d}{dw}[a_{k}]
\end{align}
and
\begin{align}
\label{14}
c_{k}=w\frac{d}{dw}[b_{k}].
\end{align}
By \eqref{13} and \eqref{14}, we obtain
\begin{align}
\label{1/1-w}
\sqrt{\frac{a_{k}c_{k}-b_{k}^{2}}{wa_{k}^{2}}} &= \sqrt{\frac{d}{dw}\left[\frac{b_{k}}{a_{k}}\right]}  \notag \\
&=\sqrt{\frac{d}{dw}\left[w\cdot \frac{d}{dw}\left[\log a_{k}\right]\right]}    \notag      \\
&=\frac{1}{(w-1)^{2}}\left(1 - \frac{(k+1)^{2}w^{k}}{a_{k}^{2}}\right).
\end{align}
Finally, we have
\begin{align}
\label{cnbm}
c_{n}b_{m}-c_{m}b_{n}&=\sum_{j=0}^{n}j^{2}w^{j}\sum_{i=0}^{m}iw^{i}-\sum_{j=0}^{m}j^{2}w^{j}\sum_{i=0}^{n}iw^{i} \notag \\
&=\sum_{j=m+1}^{n}j^{2}w^{j}\sum_{i=0}^{m}iw^{i}-\sum_{j=0}^{m}j^{2}w^{j}\sum_{i=m+1}^{n}iw^{i} \notag \\
&=\sum_{j=m+1}^{n}j^{2}w^{j}\sum_{i=0}^{m}iw^{i}-\sum_{i=0}^{m}i^{2}w^{i}\sum_{j=m+1}^{n}jw^{j} \notag \\
&=\sum_{j=m+1}^{n}\sum_{i=0}^{m}[(j-i)ijw^{i+j}]>0.
\end{align}
\subsection{Proof of \eqref{1}}
Initially, we have an estimation 
\begin{align}
\label{xiao1}
\frac{\sqrt{r_{1}r_{2}-r_{12}^{2}}}{wr_{3}^{2}}&=\sqrt{\frac{(a_{n}+a_{m})c_{n}c_{m}-c_{n}b_{m}^{2}-c_{m}b_{n}^{2}}{w^{2}(a_{n}+a_{m})^{3}}}  \notag    \\
&=\sqrt{\frac{1}{w}}\sqrt{\frac{c_{m}}{a_{n}+a_{m}}}\sqrt{\frac{(a_{n}+a_{m})c_{n}-b_{n}^{2}-\frac{c_{n}}{c_{m}}b_{m}^{2}}{w(a_{n}+a_{m})^{2}}}  \notag   \\
&\leq \sqrt{\frac{1}{w}}\sqrt{\frac{c_{m}}{a_{n}+a_{m}}}\sqrt{\frac{(a_{n}+a_{m})c_{n}-b_{n}^{2}-b_{n}b_{m}}{w(a_{n}+a_{m})^{2}}}  \notag\\
&\leq I_{1}+I_{2},
\end{align}
where $I_{1}=\sqrt{\frac{1}{w}}\sqrt{\frac{c_{m}}{a_{n}+a_{m}}}\sqrt{\frac{a_{n}c_{n}-b_{n}^{2}}{w(a_{n}+a_{m})^{2}}}$ and
$I_{2}=\sqrt{\frac{1}{w}}\sqrt{\frac{c_{m}}{a_{n}+a_{m}}}\sqrt{\frac{a_{m}c_{n}-c_{n}b_{m}}{w(a_{n}+a_{m})^{2}}}$. In the above formula, we utilize (\ref{cnbm}) in the first equality, and use the fact that $\sqrt{a+b}\leq\sqrt{a}+\sqrt{b}$, where $a$ and $b$ are positive constants, in the second equality.
Since
\begin{align*}
    \sqrt{\frac{a_{m}c_{n}-b_{n}b_{m}}{wa_{m}^{2}}} &= \sqrt{\frac{d}{dw}\left[\frac{b_{n}}{a_{m}}\right]} = \sqrt{\frac{d}{dw}\left[\frac{b_{n}a_{n}}{a_{n}a_{m}}\right]} \\
    &\leq \sqrt{\frac{a_{n}}{a_{m}}\cdot \frac{d}{dw}\left[\frac{b_{n}}{a_{n}}\right]} + \sqrt{\frac{b_{n}}{a_{n}}\cdot \frac{d}{dw}\left[\frac{a_{n}}{a_{m}}\right]},
\end{align*}
we have
\begin{align}
\label{xiao2}
    I_{2} &= \sqrt{\frac{1}{w}}\sqrt{\frac{c_{m}}{a_{n}+a_{m}}}\sqrt{\frac{a_{m}c_{n}-b_{n}b_{m}}{wa_{m}^{2}}}\sqrt{\frac{a_{m}^{2}}{(a_{n}+a_{m})^{2}}} \notag \\
    &\leq \sqrt{\frac{1}{w}}\sqrt{\frac{c_{m}}{a_{n}+a_{m}}}\sqrt{\frac{d}{dw}\left[\frac{b_{n}}{a_{n}}\right]}\sqrt{\frac{a_{n}a_{m}}{(a_{n}+a_{m})^{2}}} \notag   \\
    &\quad + \sqrt{\frac{1}{w}}\sqrt{\frac{c_{m}}{a_{n}+a_{m}}}\sqrt{\frac{d}{dw}\left[\frac{a_{n}}{a_{m}}\right]}\sqrt{\frac{b_{n}a_{m}^{2}}{(a_{n}+a_{m})^{2}}} \notag  \\
    &= \sqrt{\frac{a_{m}}{a_{n}}} \cdot I_{1} + I_{3}, 
\end{align}
where $I_{3}=\sqrt{\frac{1}{w}}\sqrt{\frac{c_{m}}{a_{n}+a_{m}}}\sqrt{\frac{d}{dw}\left[\frac{a_{n}}{a_{m}}\right]}\sqrt{\frac{b_{n}a_{m}^{2}}{(a_{n}+a_{m})^{2}}}$. We have
\begin{align*}
\sqrt{\frac{d}{dw}\left[\frac{a_{n}}{a_{m}}\right]}=\sqrt{\frac{(n-m)w^{n+m+1}-(n+1)w^{n}+(m+1)w^{m}}{(w^{m+1}-1)^{2}}}
=\sqrt{w^{m}\cdot \frac{f(w)}{g(w)}} ,
\end{align*}
where \( f(w) = (n-m)w^{n+1} - (n+1)w^{n-m} + (m+1) \) is a monotonically decreasing function with respect to \( w \), and \( g(w) = (w^{m+1}-1)^{2} \) is also monotonically decreasing. Therefore,
\begin{align}
    I_{3} & = \sqrt{w^{m-1}}\sqrt{\frac{c_{m}}{a_{n}+a_{m}}}\sqrt{\frac{f(w)}{g(w)}}\sqrt{\frac{b_{n}a_{m}^{2}}{a_{n}(a_{n}+a_{m})^{2}}} \notag \\
    & \leq \sqrt{w^{\alpha n-1}} \cdot \alpha n^{\frac{3}{2}} \cdot \sqrt{\frac{f(0)}{g\left(1-\frac{1}{n}\right)}}\notag \\
    & \leq \sqrt{w^{\alpha n-1}} \cdot \alpha n^{\frac{3}{2}} \cdot\frac{ \sqrt{\alpha n+1}}{1-e^{-\alpha}}. \notag
\end{align}
The latter inequality originates from the fact that $\left((1-\frac{1}{n})^{\alpha n+1}-1\right)^{2}\geq (e^{-\alpha}-1)^{2}$. Furthermore,
\begin{align}
\label{xiao3}
\int_{0}^{1-\frac{1}{n}}I_{3} \, dw \leq \frac{2}{\alpha n+1}e^{\frac{\alpha}{2}} \cdot \alpha n^{\frac{3}{2}} \cdot \sqrt{\frac{\alpha n+1}{e^{\alpha}-1}} \leq \frac{2e^{\frac{\alpha}{2}}\sqrt{\alpha}}{1-e^{-\alpha}}n.
\end{align}
For$(1+\sqrt{\frac{a_{m}}{a_{n}}})I_{1}$, combining \eqref{11}, \eqref{12} and  \eqref{1/1-w}, we have
\begin{align}
\label{I1}
I_{1}\leq\sqrt{\frac{1}{w}}\sqrt{\frac{c_{m}}{a_{m}}}\sqrt{\frac{a_{n}c_{n}-b_{n}^{2}}{wa_{n}^{2}}}\leq\sqrt{\frac{1}{w(1-w)^{2}}}\sqrt{\frac{c_{m}}{a_{m}}}.
\end{align}
Over the interval $0\leq w\leq \frac{1}{2}$, the following estimate can be obtained,
\[I_{1}\leq 2\alpha n\sqrt{\frac{1}{w}},\]
which implies
\begin{align}
\label{xiao4}
\int_{0}^{\frac{1}{2}}\left(1+\sqrt{\frac{a_{m}}{a_{n}}}\right)I_{1}dw\leq \int_{0}^{\frac{1}{2}}2\alpha n\sqrt{\frac{1}{w}}dw=4\sqrt{2}\cdot \alpha n.
\end{align}
For the remaining interval $\frac{1}{2} \leq w \leq 1-\frac{1}{n}$, we make a more precise estimate of $\sqrt{\frac{c_{k}}{a_{k}}}$. We begin by performing polynomial manipulation on $c_{k}$ to achieve the desired form,
\begin{align*}
c_{k}-wc_{k}=\sum_{j=0}^{k}j^{2}w^{j}-w\sum_{j=0}^{k}j^{2}w^{j}=\sum_{j=1}^{k}(2j-1)w^{j}-k^{2}w^{k+1}.
\end{align*}
After another similar operation, we have
\begin{align*}
&\quad(1-w)^{2}c_{k}+k^{2}(w^{k+1}-w^{k+2})\\
&=c_{k}-wc_{k}+k^{2}w^{k+1}-k(c_{k}-wc_{k}+k^{2}w^{k+1})\\
&=\sum_{j=1}^{k}(2j-1)w^{j}-w\sum_{j=1}^{k}(2j-1)w^{j}\\
&=\sum_{j=1}^{k}2w^{j}-w-(2k-1)w^{k+1}.
\end{align*}
Then divide $c_{k}$ by $a_{k}$,
\begin{align}
\label{xiao5}
\frac{c_{k}}{a_{k}} &=\left(\frac{\sum_{j=1}^{k}2w^{j}-w-(2k-1)w^{k+1}-k^{2}(w^{k+1}-w^{k+2})}{(1-w)^{2}}\right)\Big/\left(\frac{w^{k+1}-1}{w-1}\right)\notag\\
&= \frac{w(k^{2}w^{k}(w-1)^{2}+2k(w^{k}-w^{k+1})+(w^{k}-1)(w+1))}{(w-1)^{2}(w^{k+1}-1)} \notag\\
&\leq \frac{w(w+1)}{(1-w)^{2}} - \frac{2kw^{k+1}}{(1-w)(1-w^{k+1})} - \frac{k^{2}w^{k+1}}{1-w^{k+1}}\notag\\
&=\frac{2}{(1-w)^2}-\frac{3}{1-w}+1- \frac{2kw^{k+1}}{(1-w)(1-w^{k+1})} - \frac{k^{2}w^{k+1}}{1-w^{k+1}}. \notag
\end{align}
We utilize the fact that $\frac{w^{k}-1}{w^{k+1}-1}<1$ when $w<1$ in the above inequality. Clearly, for all terms in the above equation except $\frac{2}{(1-x)^2}$, the sum is negative when $w<1$. Thus, incorporating \eqref{I1}, we obtain the estimation
\begin{align}
\int_{\frac{1}{2}}^{1-\frac{1}{n}}\left(1+\sqrt{\frac{a_{m}}{a_{n}}}\right)I_{1}dw&\leq 4\int_{\frac{1}{2}}^{1-\frac{1}{n}}\frac{1}{(1-w)^{2}}dw=4(n-2).
\end{align}
Therefore, combining \eqref{xiao1}-\eqref{xiao3} and \eqref{xiao4}-\eqref{xiao5}, we have
\begin{align}
\label{eq:1in1}
\int_{0}^{1-\frac{1}{n}}\frac{\sqrt{r_{1}r_{2}-r_{12}^{2}}}{wr_{3}^{2}}dw\leq \frac{2e^{\frac{\alpha}{2}}\sqrt{\alpha}}{1-e^{-\alpha}}n+4\sqrt{2}\alpha n+4(n-2),
\end{align}
Meanwhile, by \eqref{11} we have                                                                                                                                                                                                                                                                                                                                                                                                                                                                                                                                                                                                             
\begin{align}
\label{eq:2in1}
\int_{0}^{1-\frac{1}{n}}\frac{r_{1}-r_{2}}{wr_{3}^{2}} \, dw & = \int_{0}^{1-\frac{1}{n}}\frac{(a_{n}+a_{m})(c_{n}-c_{m})-(b_{n}^{2}-b_{m}^{2})}{w(a_{n}+a_{m})^{2}} \, dw \notag\\
& = \frac{b_{n}-b_{m}}{a_{n}+a_{m}}\bigg\vert_{0}^{1-\frac{1}{n}} \notag\\
& = \frac{\sum_{j=0}^{n}j(1-\frac{1}{n})^{j}-\sum_{j=0}^{m}(1-\frac{1}{n})^{j}}{\sum_{j=0}^{n}(1-\frac{1}{n})^{j}+\sum_{j=0}^{m}(1-\frac{1}{n})^{j}}\notag \\
& \leq\frac{\sum_{j=0}^{n}j(1-\frac{1}{n})^{j}}{\sum_{j=0}^{n}(1-\frac{1}{n})^{j}}\leq n.
\end{align}
Combining \eqref{F}, $(\ref{eq:1in1})$ and $(\ref{eq:2in1})$, we complete the proof of (\ref{1}).
\subsection{Proof of \eqref{2}}
We first make a transformation of $F(w)$,
\begin{align*}
F(w)= \frac{r_{1}^{2}+r_{2}^{2}-2r_{12}^{2}}{wr_{3}^{2}\sqrt{(r_{1}+r_{2})^{2}-4r_{12}^{2}}}=\frac{\sqrt{(r_{1}+r_{2})^{2}-4r_{12}^{2}}}{wr_{3}^{2}}-\frac{2r_{1}r_{2}-2r_{12}^{2}}{wr_{3}^{2}\sqrt{(r_{1}+r_{2})^{2}-4r_{12}^{2}}}.
\end{align*}
Since
\begin{align*}
r_{1}r_{2}-r_{12}^{2} &= \left(a_{n}+a_{m}\right)^{2}c_{n}c_{m} - \left(a_{n}+a_{m}\right)\left(c_{n}b_{m}^{2} + c_{m}b_{n}^{2}\right) \\
&= \left(a_{n}+a_{m}\right)\left[c_{n}\left(a_{m}c_{m}-b_{m}^{2}\right)+c_{m}\left(a_{n}c_{n}-b_{n}^{2}\right)\right] \\
&= \left(\sum_{j=0}^{n}w^{j}+\sum_{j=0}^{m}w^{j}\right)\left[\sum_{j=0}^{n}j^{2}w^{j}\left(\sum_{j=0}^{m}w^{j}\sum_{j=0}^{m}j^{2}w^{j}-\sum_{j=0}^{m}jw^{j}\sum_{j=0}^{m}jw^{j}\right)\right.\\
&\quad + \left.\sum_{j=0}^{m}j^{2}w^{j}\left(\sum_{j=0}^{n}w^{j}\sum_{j=0}^{n}j^{2}w^{j}-\sum_{j=0}^{n}jw^{j}\sum_{j=0}^{n}jw^{j}\right)\right] \\
&= \left(\sum_{j=0}^{n}w^{j}+\sum_{j=0}^{m}w^{j}\right)\left[\sum_{j=0}^{n}j^{2}w^{j}\sum_{j=0}^{m}\sum_{i=0}^{m}\left(\frac{i^{2}+j^{2}}{2}-ij\right)w^{i+j}\right.\\
&\quad + \left.\sum_{j=0}^{m}j^{2}w^{j}\sum_{j=0}^{n}\sum_{i=0}^{n}\left(\frac{i^{2}+j^{2}}{2}-ij\right)w^{i+j}\right] \\
&\geq 0,
\end{align*}
we have
\begin{align}
\label{F2}
F(w)\leq \frac{\sqrt{(r_{1}+r_{2})^{2}-4r_{12}^{2}}}{wr_{3}^{2}}\leq\frac{r_{1}+r_{2}}{wr_{3}^{2}}=I_{4}+I_{5},
\end{align}
where $I_{4}=\frac{r_{1}+r_{2}-2r_{12}}{wr_{3}^{2}}$ and $I_{5}=\frac{2r_{12}}{wr_{3}^{2}}$.
Start by calculating the integral of $I_{4}$,
\begin{align}
\label{I4}
\int_{1-\frac{1}{n}}^{1+\frac{1}{n}} I_{4} \, dw 
&= \int_{1-\frac{1}{n}}^{1+\frac{1}{n}} \frac{(a_{n}+a_{m})(c_{n}+c_{m}) - (b_{n}+b_{m})^{2}}{w(a_{n}+a_{m})^{2}} \, dw \notag \\
&= \frac{b_{n}+b_{m}}{a_{n}+a_{m}} \bigg \vert_{1-\frac{1}{n}}^{1+\frac{1}{n}} \notag \\
&= \frac{\sum_{j=0}^{n}jw^{j} + \sum_{j=0}^{m}jw^{j}}{\sum_{j=0}^{n}w^{j} + \sum_{j=0}^{m}w^{j}} \bigg \vert_{1-\frac{1}{n}}^{1+\frac{1}{n}} \notag \\
&\leq \frac{\sum_{j=0}^{n}j\left(1+\frac{1}{n}\right)^{j} + \sum_{j=0}^{m}j\left(1+\frac{1}{n}\right)^{j}}{\sum_{j=0}^{n}\left(1+\frac{1}{n}\right)^{j} + \sum_{j=0}^{m}\left(1+\frac{1}{n}\right)^{j}} \notag \\
&\leq \frac{\sum_{j=0}^{n}j\left(1+\frac{1}{n}\right)^{j}}{\sum_{j=0}^{n}\left(1+\frac{1}{n}\right)^{j}} + \frac{\sum_{j=0}^{m}j\left(1+\frac{1}{n}\right)^{j}}{\sum_{j=0}^{m}\left(1+\frac{1}{n}\right)^{j}} \notag \\
&\leq (1+\alpha) n. 
\end{align}
In the last inequality of the above formula, we use \eqref{11}. Then we calculate the integral of $I_{5}$,
\begin{align}
\label{I5}
\int_{1-\frac{1}{n}}^{1+\frac{1}{n}}I_{5}dw&=2\int_{1-\frac{1}{n}}^{1+\frac{1}{n}}\frac{b_{n}b_{m}}{w(a_{n}+a_{m})^{2}}dw\notag\\
&\leq2\int_{1-\frac{1}{n}}^{1+\frac{1}{n}}\frac{b_{n}b_{m}}{4wa_{n}a_{m}}dw\notag\\
&\leq\frac{\alpha n^{2}}{2}\int_{1-\frac{1}{n}}^{1+\frac{1}{n}}\frac{1}{w}dw \notag\\
&=\frac{\alpha n^{2}}{2}\log \frac{n+1}{n-1}\notag \\
&=\frac{\alpha n}{2}\log \left(1+\frac{2}{n-1}\right)^{n} \notag\\
&\leq \alpha n.
\end{align}
Combining \eqref{F2}-\eqref{I5}, we complete the proof of \eqref{2}.
\subsection{Proof of \eqref{3}}
We first make an estimation as follows,
\begin{align}
\label{41}
\frac{\sqrt{r_{1}r_{2}-r_{12}^{2}}}{wr_{3}^{2}} & = \sqrt{\frac{(a_{n}+a_{m})c_{n}c_{m}-c_{n}b_{m}^{2}-c_{m}b_{n}^{2}}{w^{2}(a_{n}+a_{m})^{3}}}\notag \\
&\leq \sqrt{\frac{a_{n}c_{n}c_{m}-c_{m}b_{n}^{2}}{w^{2}(a_{n}+a_{m})^{3}}}+\sqrt{\frac{a_{m}c_{n}c_{m}-c_{n}b_{m}^{2}}{w^{2}(a_{n}+a_{m})^{3}}}\notag \\
&=\sqrt{\frac{1}{w}}\sqrt{\frac{c_{m}}{a_{m}}}\sqrt{\frac{a_{n}c_{n}-b_{n}^{2}}{wa_{n}^{2}}}\sqrt{\frac{a_{n}^{2}a_{m}}{(a_{n}+a_{m})^{3}}}\notag \\
&\quad +\sqrt{\frac{1}{w}}\sqrt{\frac{c_{n}}{a_{n}}}\sqrt{\frac{a_{m}c_{m}-b_{m}^{2}}{wa_{m}^{2}}}\sqrt{\frac{a_{n}a_{m}^{2}}{(a_{n}+a_{m})^{3}}}\notag \\
&\leq \sqrt{\frac{1}{w}}\left( \frac{\alpha n}{w-1}\sqrt{\frac{a_{m}}{a_{n}+a_{m}}}+\frac{n}{w-1}\sqrt{\frac{a_{m}}{a_{n}+a_{m}}}\right)\notag \\
&=\frac{(1+\alpha)n}{\sqrt{w}(w-1)}\sqrt{\frac{a_{m}}{a_{n}+a_{m}}}.
\end{align}
We utilize \eqref{12} and \eqref{1/1-w} in the above second inequality.
When $w>1+\frac{\log n}{n}$, by the monotonicity of the function, we establish 
\begin{align}
\label{4x}
\frac{a_{m}}{a_{n}+a_{m}}&<\frac{a_{m}}{a_{n}}=\frac{w^{m}-1}{w^{n}-1}\leq w^{(\alpha-1)n}\leq 2 n^{\alpha-1}.
\end{align}
The last inequality above holds for sufficiently large $n$. Combining \eqref{41} and \eqref{4x}, we have\\
\begin{align}
\label{42}
\int_{1+\frac{\log n}{n}}^{2}\frac{\sqrt{r_{1}r_{2}-r_{12}^{2}}}{wr_{3}^{2}} \, dw & \leq \int_{1+\frac{\log n}{n}}^{2}\frac{(1+\alpha)\sqrt{n}}{\sqrt{w}(w-1)}\sqrt{\frac{a_{m}}{a_{n}+a_{m}}} \, dw\notag \\
& \leq \int_{1+\frac{\log n}{n}}^{2}2\frac{(1+\alpha)n^{\frac{1+\alpha}{2}}}{\sqrt{w}(w-1)} \, dw\notag \\
& \leq \int_{1+\frac{\log n}{n}}^{2}2\frac{(1+\alpha)n^{\frac{1+\alpha}{2}}}{(w-1)} \, dw\notag \\
& = 2(1+\alpha)n^{\frac{1+\alpha}{2}}\log(w-1) \bigg\vert_{1+\frac{\log n}{n}}^{2}\notag \\
& = 2(1+\alpha)n^{\frac{1+\alpha}{2}}\log\left(\frac{n}{\log n}\right)\notag \\
& <n.
\end{align}
Similarly, for $w>2$,
\begin{align}
\label{43}
\int_{2}^{\infty}\frac{\sqrt{r_{1}r_{2}-r_{12}^{2}}}{wr_{3}^{2}} \, dw&\leq\int_{2}^{\infty}\frac{(1+\alpha)n^{\frac{1+\alpha}{2}}}{\sqrt{w}(w-1)} \, dw\notag \\
&\leq\int_{2}^{\infty}2\frac{(1+\alpha)n^{\frac{1+\alpha}{2}}}{(w-1)^{\frac{3}{2}}} \, dw\notag \\
&\leq 4(1+\alpha)n^{\frac{1+\alpha}{2}}\notag \\
&<n.
\end{align}
Over the remaining interval $1+\frac{1}{n}\leq w\leq 1+\dfrac{\log n}{n}$, we have a more precise estimate of $\sqrt{\frac{a_{m}}{a_{m}+a_{n}}}$,
\begin{align}
\label{4x2}
\sqrt{\frac{a_{m}}{a_{m}+a_{n}}}=\sqrt{\frac{1}{\frac{w^{n+1}-1}{w^{m+1}-1}+1}}\leq\sqrt{\frac{1}{w^{(1-\alpha)n}+1}}\leq \sqrt{\frac{1}{w^{(1-\alpha)n}}}.
\end{align}
Therefore, plugging \eqref{4x2} into \eqref{41}, we have
\begin{align*}
\int_{1+\frac{1}{n}}^{1+\frac{\log n}{n}} \frac{\sqrt{r_{1}r_{2}-r_{12}^{2}}}{wr_{3}^{2}} \, dw  \leq (1+\alpha)n \int_{1+\frac{1}{n}}^{1+\frac{\log n}{n}} \frac{1}{w-1} \sqrt{\frac{1}{w^{(1-\alpha)n}}} \, dw.\\
\end{align*}
Replacing $w$ with $1+\dfrac{t}{n}$, where $1\leq t\leq \log n$, yields
\begin{align}
\label{44}
\int_{1+\frac{1}{n}}^{1+\frac{\log n}{n}} \frac{\sqrt{r_{1}r_{2}-r_{12}^{2}}}{wr_{3}^{2}} \, dw 
&\leq (1+\alpha)n \int_{1}^{\log n} \frac{1}{t} \sqrt{\frac{1}{\left(1+\frac{t}{n}\right)^{(1-\alpha)n}}} \, dt\notag \\
&\leq (1+\alpha)n \int_{1}^{\log n} \frac{1}{t} \sqrt{\frac{1}{\exp \left((1-\alpha)\left(t-\frac{t^{2}}{2n}\right)\right)}} \, dt\notag \\
&\leq (1+\alpha)n \int_{1}^{\log n} \frac{1}{t^{\frac{3-\alpha}{2}}} \, dt\notag \\
&=\frac{2(1+\alpha)n}{1-\alpha}\left(1-\left(\log n\right)^{\frac{\alpha-1}{2}}\right)\notag \\
&\leq \frac{2(1+\alpha)n}{1-\alpha},
\end{align}
where we utilize the fact that $\exp\left( t - \frac{t^2}{2n} \right) \geq t$ when $1\leq t\leq \log n$ in the third inequality. Thus, combining \eqref{42}, \eqref{43} and \eqref{44}, we have 
\begin{align}
\label{45}
\int_{1+\frac{1}{n}}^{\infty} \frac{\sqrt{r_{1}r_{2}-r_{12}^{2}}}{wr_{3}^{2}} \, dw &\leq \left[2+\frac{2(1+\alpha)}{1-\alpha}\right]n.
\end{align}
Meanwhile, we have
\begin{align}
\label{46}
\int_{1+\frac{1}{n}}^{\infty}\frac{r_{1}-r_{2}}{wr_{3}^{2}} \, dw & = \int_{1+\frac{1}{n}}^{\infty}\frac{(a_{n}+a_{m})(c_{n}-c_{m})-(b_{n}^{2}-b_{m}^{2})}{w(a_{n}+a_{m})^{2}} \, dw \notag\\
& = \frac{b_{n}-b_{m}}{a_{n}+a_{m}}\bigg\vert_{1+\frac{1}{n}}^{\infty} \notag\\
& \leq \lim_{w\rightarrow\infty}\frac{b_{n}-b_{m}}{a_{n}+a_{m}} \notag \\
&=n.
\end{align}
Combining \eqref{F}, \eqref{45} and \eqref{46}, we complete the proof of \eqref{3}.

\subsection{Proof of \eqref{4}}
Finally, we make an estimate of the lower bound.
Since
\begin{align*}
F(w)&=\frac{r_{1}^{2}+r_{2}^{2}-2r_{12}^{2}}{wr_{3}^{2}\sqrt{(r_{1}+r_{2})^{2}-4r_{12}^{2}}} \geq\frac{r_{1}-r_{2}}{wr_{3}^{2}} =\frac{d}{dw}\left[\frac{b_{n}-b_{m}}{a_{n}+a_{m}}\right],
\end{align*}
we have
\begin{align}
\label{51}
\int_{0}^{\infty}F(w)dw &\geq\int_{0}^{\infty}\frac{r_{1}-r_{2}}{wr_{3}^{2}}dw=\frac{b_{n}-b_{m}}{a_{n}+a_{m}}\bigg\vert_{0}^{\infty}=\lim_{w\rightarrow\infty}\frac{b_{n}-b_{m}}{a_{n}+a_{m}} = n.
\end{align}
Combining \eqref{1}-\eqref{3} and \eqref{51}, we have
\[ n\leq \mathbb{E}[N_{H}(\mathbb{C})]\leq C_{0}n, \]
where $C_{0}$ is a constant related to $\alpha$. Hence we finish the proofs.

\bibliography{reference}

\providecommand{\bysame}{\leavevmode\hbox to3em{\hrulefill}\thinspace}
\providecommand{\MR}{\relax\ifhmode\unskip\space\fi MR }
\providecommand{\MRhref}[2]{%
  \href{http://www.ams.org/mathscinet-getitem?mr=#1}{#2}
}
\providecommand{\href}[2]{#2}
\begin{thebibliography}{10}

\bibitem{bleher2014counting}
P.~M. Bleher, Y.~Homma, L.~L. Ji, and R.~K.~W. Roeder, \emph{Counting zeros of
  harmonic rational functions and its application to gravitational lensing},
  International Mathematics Research Notices \textbf{2014} (2014), no.~8,
  2245--2264.

\bibitem{Bloch1932}
A.~Bloch and G.~Pólya, \emph{On the roots of certain algebraic equations},
  Proceedings of the London Mathematical Society \textbf{s2-33} (1932),
  102--114.

\bibitem{bshouty1995exact}
D.~Bshouty, W.~Hengartner, and T.~Suez, \emph{The exact bound on the number of
  zeros of harmonic polynomials}, Journal d’Analyse Math{\'e}matique
  \textbf{67} (1995), 207--218.

\bibitem{bshouty2004crofoot}
D.~Bshouty and A.~Lyzzaik, \emph{On crofoot-sarason's conjecture for harmonic
  polynomials}, Computational Methods and Function Theory \textbf{4} (2004),
  35--42.

\bibitem{dunnage1966number}
J.~Dunnage, \emph{The number of real zeros of a random trigonometric
  polynomial}, Proceedings of the London Mathematical Society \textbf{3}
  (1966), 53--84.

\bibitem{edelman1995many}
A.~Edelman and E.~Kostlan, \emph{How many zeros of a random polynomial are
  real?}, Bulletin of the American Mathematical Society \textbf{32} (1995),
  1--37.

\bibitem{geyer2008sharp}
L.~Geyer, \emph{Sharp bounds for the valence of certain harmonic polynomials},
  Proceedings of the American Mathematical Society \textbf{136} (2008),
  549--555.

\bibitem{ibragimov1997roots}
I.~Ibragimov and O.~Zeitouni, \emph{On roots of random polynomials},
  Transactions of the American Mathematical Society \textbf{349} (1997),
  2427--2441.

\bibitem{kac1943average}
M.~Kac, \emph{On the average number of real roots of a random algebraic
  equation}, Bulletin of the American Mathematical Society \textbf{49} (1943),
  314--320.

\bibitem{khavinson2010transcendental}
D.~Khavinson and E.~Lundberg, \emph{Transcendental harmonic mappings and
  gravitational lensing by isothermal galaxies}, Complex Analysis and Operator
  Theory \textbf{4} (2010), no.~3, 515--524.

\bibitem{khavinson2024valence}
D.~Khavinson, E.~Lundberg, and S.~Perry, \emph{On the valence of logharmonic
  polynomials}, Recent Progress in Function Theory and Operator Theory
  \textbf{799} (2024), 23--40.

\bibitem{khavinson2006number}
D.~Khavinson and G.~Neumann, \emph{On the number of zeros of certain rational
  harmonic functions}, Proceedings of the American Mathematical Society
  \textbf{134} (2006), 1077--1085.

\bibitem{Khavinson2002}
D.~Khavinson and G.~Światek, \emph{On the number of zeros of certain harmonic
  polynomials}, Proceedings of the American Mathematical Society \textbf{131}
  (2003), 409--414.

\bibitem{lee2017new}
S.~Lee and A.~Saez, \emph{A new lower bound for the maximal valence of harmonic
  polynomials}, Computational Methods and Function Theory \textbf{17} (2017),
  139--149.

\bibitem{Lerario2016}
A.~Lerario and E.~Lundberg, \emph{On the zeros of random harmonic polynomials:
  the truncated model}, Journal of Mathematical Analysis and Applications
  \textbf{438} (2016), 1041--1054.

\bibitem{Li2009}
W.~V. Li and A.~Wei, \emph{On the expected number of zeros of a random harmonic
  polynomial}, Proceedings of the American Mathematical Society \textbf{137}
  (2009), 195--204.

\bibitem{lundberg2023valence}
E.~Lundberg, \emph{The valence of harmonic polynomials viewed through the
  probabilistic lens}, Proceedings of the American Mathematical Society
  \textbf{151} (2023), 2963--2973.

\bibitem{Lundberg2023}
E.~Lundberg and A.~Thomack, \emph{On the average number of zeros of random
  harmonic polynomials with iid coefficients: precise asymptotics}, arXiv
  preprint arXiv:2308.10333 (2023).

\bibitem{Peres2005}
Y.~Peres and B.~Virág, \emph{Zeros of the i.i.d. gaussian power series: a
  conformally invariant determinantal process}, Acta Mathematica \textbf{194}
  (2005), 1--35.

\bibitem{petters2009mathematical}
A.~O. Petters, B.~Rider, and A.~M. Teguia, \emph{{A mathematical theory of
  stochastic microlensing. I. Random time delay functions and lensing maps}},
  Journal of Mathematical Physics \textbf{50} (2009), no.~7, 072503.

\bibitem{petters2009mathematical2}
\bysame, \emph{{A mathematical theory of stochastic microlensing. II. Random
  images, shear, and the Kac–Rice formula}}, Journal of Mathematical Physics
  \textbf{50} (2009), no.~12, 122501.

\bibitem{shepp1995complex}
L.~A. Shepp and R.~J. Vanderbei, \emph{The complex zeros of random
  polynomials}, Transactions of the American Mathematical Society (1995),
  4365--4384.

\bibitem{Thomack2016}
A.~Thomack, \emph{On the zeros of random harmonic polynomials: the naive
  model}, arXiv preprint arXiv:1610.02611 (2016).

\bibitem{Thomack2018}
A.~Thomack and Z.~Tyree, \emph{On the zeros of random harmonic polynomials: the
  weyl model}, Analysis and Mathematical Physics \textbf{8} (2018), 237--253.

\bibitem{wei2017numbers}
A.~Wei, \emph{{On the numbers of images of two stochastic gravitational lensing
  models}}, Journal of Mathematical Physics \textbf{58} (2017), no.~2, 022501.

\bibitem{Wilmshurst1994}
A.~S. Wilmshurst, \emph{Complex harmonic mappings and the valence of harmonic
  polynomials}, D. phil. thesis, University of York, England, 1994.

\bibitem{wilmshurst1998valence}
\bysame, \emph{The valence of harmonic polynomials}, Proceedings of the
  American Mathematical Society \textbf{126} (1998), 2077--2081.

\end{thebibliography}

\end{document}